\documentclass{article}
\advance\voffset by -25mm
\advance\hoffset by -25mm
\usepackage[latin1]{inputenc}
\usepackage{amsmath}
\usepackage{amsfonts}
\usepackage{verbatim}

\newtheorem{prop}{Proposition}

\title{}
\begin{document}

\title
{Expanding graphs, Ramanujan graphs, and $1$-factor perturbations}

\author{
Pierre de la Harpe
and
Antoine Musitelli\\adress email:  Pierre.delaHarpe@math.unige.ch 
and musitel0@math.unige.ch}
\date{\today}
\maketitle

\abstract
We construct $(k \pm 1)$-regular graphs
which provide sequences of expanders
by adding or substracting appropriate $1$-factors
from given sequences of $k$-regular graphs.
We compute numerical examples in a few cases for which
the given sequences are from the work of Lubotzky, Phillips, and Sarnak
(with $k-1$ the order of a finite field).
If $k+1 = 7$, our construction results in a sequence of $7$-regular
expanders with all spectral gaps at least $6 - 2\sqrt 5 \approx 1.52$;
the corresponding minoration for a sequence of
Ramanujan $7$-regular graphs (which is not known to exist)
would be $7 - 2\sqrt 6 \approx 2.10$.
\endabstract

\section{Introduction}

  Let $X = (V,E)$ be a simple finite graph with $n$ vertices, where
$V$ denotes the vertex set and $E$ the set of geometrical edges of $X$.
The adjacency matrix $A$ of $X$, with rows and columns indexed by $V$,
is defined by $A_{v,w} = 1$ if there exists an edge connecting $v$ and
$w$, and $A_{v,w} = 0$ otherwise (in particular $A_{v,v} = 0$).
The eigenvalues of $X$, which are those of $A$,
constitute a decreasing sequence
$\lambda_0(X) \ge \lambda_1(X) \ge  \hdots \ge \lambda_{n-1}(X)$;
the {\it spectral gap} $\lambda_0(X) - \lambda_1(X)$ of $X$ is positive
if and only if $X$ is connected.
Let us assume from now on that $X$ is $k$-regular for some $k \ge 3$,
namely that $\sum_{w}A_{v,w} = k$ for all $v \in V$,
so that $\lambda_0(X) = k$.

Recall that, for any infinite sequence $(X_i)_{i \in I}$ of
connected $k$-regular simple finite graphs with increasing vertex sizes,
we have $\liminf_{i \to \infty} \lambda_1(X_i) \ge 2\sqrt{k-1}$.
A graph $X$ is said to be a {\it Ramanujan graph} if
it is connected and if $|\mu| \le 2\sqrt{k-1}$
for any eigenvalue $\mu \ne \pm k$ of $X$.
 From elaborate arithmetic constructions,
we know explicit infinite sequences of Ramanujan graphs
for degree $k$ when $k-1$ is the order of a finite field;
but the existence of such sequences is an open problem
for other degrees, for example when $k = 7$.
It is thus interesting to find {\it sequences of expanders of degree
$k$},  namely infinite sequences
$(X_i)_{i \in I}$ of $k$-regular connected simple finite graphs with
increasing vertex sizes such that $\inf_{i \in I}(k-\lambda_1(X_{i}))$ is
strictly positive, and indeed as large as possible (short of being equal
to  $k-2\sqrt{k-1}$).

  For all this, see for example [Lubot--94], [Valet--97],
[Colin--98], and [DaSaV--03].

\medskip

   The object of the present Note is to examine a procedure of
construction of sequences of expanders $(X_i)_{i \in I}$
of degree $k$ by perturbation of sequences of Ramanujan graphs.
When $k-l$ is the order of a finite field, we obtain estimates
$\lambda_1(X_i) \le l-1 + 2\sqrt{(k-l)}$;  for example, for $k = 7$ and $l=1$,
this corresponds to a spectral gap
$$
7 - \lambda_1(X_i) \, \ge \,
6 - 2 \sqrt{5} \, \approx \, 1.52
\qquad \text{for all $i \in I$,}
$$
to be compared with the Ramanujan minoration by
$$
7 - \liminf_{i \in I}\lambda_1(X_i) \, \le \,
7 - 2\sqrt{6} \, \approx \, 2.10
$$
of the spectral gap.

We insist on finding explicit constructions, but we record however the 
following results of J. Friedman (see [Fried--91], and later work) based on 
random techniques: for all $k \ge 3$ and all $\epsilon > 0$, there exists 
sequences $\left( X_i \right)_{i \in I}$ of connected $k$-regular simple 
finite graphs with increasing vertex sizes and with $\lambda_1(X_i) \le 
2\sqrt{k-1}+\epsilon$ for all $i \in I$.
 
\medskip
   Let $X = (V,E)$ be a graph.
If $X$ is not bipartite,
we denote by $\overline{X} = (V,\overline E)$ the
{\it complement} of $X$;
two distinct vertices are adjacent in $\overline{X}$
if and only if they are not so in $X$.
If $X$ is bipartite, given with a bipartition $V = V_0 \sqcup V_1$,
we denote by $\overline{X} = (V,\overline E)$ the
{\it bipartite complement} of $X$;
two vertices $v \in V_0$, $w \in V_1$ are adjacent in $\overline{X}$
if and only if they are not in $X$.
A {\it matching} of  a graph $X$
is a subset $M$ of $E$ such that any vertex $x \in V$
is incident with at most one edge of $M$,
and a {\it perfect matching} (also called {\it $1$-factor})
is a subset $F$ of $E$ such that any
vertex $x \in V$ is incident with exactly one edge of $F$.

   Let $X = (V,E)$ be a graph. If $F$ is a perfect matching of $X$,
we denote by $X-F$ the graph $(V,E \setminus F)$; if $X$ is $k$-regular,
then $X-F$ is $(k-1)$-regular. If $F$ is a perfect matching of
$\overline{X}$, we denote by $X+F$ the graph
$\overline{\overline{X}-F}$; if $X$ is $k$-regular, then $X+F$ is
$(k+1)$-regular.

   The basic observation for the present Note is the set of
inequalities
$$
| \lambda_j(X \pm F) - \lambda_j(X)| \, \le \, 1
$$
for any perfect matching $F$ of $X$ (for $X-F$) or of $\overline{X}$
(for $X+F$), and for all ${j \in \{0,\hdots,n-1\}}$,
where $n = |V|$
(Proposition 2).
We can apply this to the Ramanujan graphs
$X^{p,q}$ and their complements
(notation of [DaSaV-03], see below).
In Section 3, we describe an algorithm for finding perfect matchings
in regular bipartite graphs
(thus concentrating on pairs $(p,q)$
for which the graph $X^{p,q}$ is bipartite).
In conclusion, we report some
numerical computations.

\bigskip

\section{ Graphs of the form $X^{p,q} \pm F$}

\medskip

 Let us recall the definition of the graphs $X^{p,q}$.

   If $R$ is a commutative ring with unit,
the {\it Hamilton quaternion algebra $\Bbb H (R)$ over $R$}
is the free module $R^4$ with basis $\{1,i,j,k\}$,
where multiplication is defined by
$i^2 = j^2 = k^2 = -1$, and $ij = -ji = k$,
plus circular permutations of $i,j,k$.
A quaternion $q = a_0 + a_1i + a_2j + a_3k$ has
a {\it conjugate} $\overline q = a_0 - a_1i - a_2j - a_3k$ and
a {\it norm} $N(q) = \overline q q = a_0^2 + a_1^2 + a_2^2 + a_3^2$.

   Let $p \in \Bbb N$ be an odd prime.
If $p \equiv 1 \pmod 4$,
a theorem of Jacobi shows that there are
exactly $p+1$ quaternions in $\Bbb H(\Bbb Z)$ of norm $p$
of the form $a_0 + a_1i + a_2j + a_3k$ with $a_0 \equiv 1 \pmod 2$,
and $a_0 \ge 1$.
These occur in pairs $(\alpha,\overline{\alpha})$;
we select arbitrarily one, say $\alpha_l$, from each pair,
and we set
$$
S_p \, = \, \{ \alpha_1, \overline{\alpha_1}, \hdots,
\alpha_{s}, \overline{\alpha_{s}} \}
\qquad \text{with} \quad 2s = p+1.
$$
If $p \equiv 3 \pmod 4$, there are quaternions
in $\Bbb H(\Bbb Z)$ of norm $p$
of the form $a_0 + a_1i + a_2j + a_3k$ with $a_0 \equiv 0 \pmod 2$,
and $a_0 \ge 0$.
 From those with $a_0 \ge 2$, say $2s$ of them,
we obtain $\alpha_1, \hdots, \alpha_s$ as above.
Those of the form $a_1 i + a_2 j + a_3 k$,
say $2t$ of them
\footnote{
Observe that $2t$ is a multiple of $8$,
since each of $a_1,$ $a_2$, $a_3$ is odd,
in particular not $0$,
so that each sign change provides another writing of $p$
as a sum of three squares.
} \hskip-.1cm
,
occur in pairs $(\beta,-\beta)$;
we select arbitrarily one, say $\beta_m$, from each pair,
and we set
$$
S_p \, = \, \{ \alpha_1, \overline{\alpha_1}, \hdots,
\alpha_{s}, \overline{\alpha_{s}}, \beta_1, \hdots, \beta_t \} .
$$
Observe that $t/4$ is the number of solutions in $\Bbb N$
of the equation $a_1^2 + a_2^2 + a_3^2 = p$,
and that we have again $|S_p| = 2s + t = p+1$
by Jacobi's theorem.
Observe also that we can have $s = 0$ (case of $p = 3$),
as well as $t = 0$ (case of $p \equiv 7 \pmod{8}$),
or both $s$ and $t$ positive (case of $p = 19$,
with $s = 4$ and $t = 12$).

   Let $q$ be another odd prime, $q \ne p$, and let
$\tau_q : \Bbb H (\Bbb Z) \longrightarrow \Bbb H (\Bbb F_q)$
denote reduction modulo $q$.
The equation $x^2 + y^2 + 1 = 0$ has solutions in $\Bbb F_q$.
We choose one solution; then the mapping
$\psi_q : \Bbb H (\Bbb F_q) \longrightarrow M_2(\Bbb F_q)$ defined by
$$
\psi_q(a_0 + a_1 i + a_2 j + a_3 k) \, = \,
\left(
\begin{array}{cc}
a_0 + a_1 x + a_3 y & -a_2 y + a_1 + a_3 x \\
-a_1 y - a_2 + a_3 x & a_0 - a_1 x - a_2 y
\end{array}
 \right)
$$
is an algebra isomorphism and
$\psi_q\left( \tau_q (S_p) \right)$ is in the group $GL_2(q)$
of invertible elements of $M_2(\Bbb F_q)$.
We denote by $\phi : GL_2(q) \longrightarrow PGL_2(q)$
the reduction modulo the centre, and we set
$$
S_{p,q} \, = \, \phi
\Big( \psi_q \left( \tau_q ( S_p ) \right) \Big)
\, \subset \, PGL_2(q) .
$$
It follows from the definitions that $S_{p,q}$ is symmetric:
if $s \in S_{p,q}$ is the image of $\alpha_l \in S_p$
(notation as above), then $s^{-1}$ is the image of
$\overline{\alpha_l}$;
if $s$ is the image of $\beta_m \in S_p$, then $s^2 = 1$.
Moreover, it is known that $|S_{p,q}| = p+1$.
There are now two cases to consider.

   Either $p$ is a square modulo $q$.
Then $S_{p,q} \subset PSL_2(q)$
and indeed $S_{p,q}$ generates $PSL_2(q)$.
By definition, $X^{p,q}$ is the Cayley graph
of $PSL_2(q)$ with respect to $S_{p,q}$;
more precisely, $X^{p,q} = (V,E)$ with $V = PSL_2(q)$
and $\{v,w\} \in E$ if $v^{-1}w  \in S^{p,q}$.
It is a $(p+1)$-regular graph
with $\frac{1}{2}q(q^2-1)$ vertices
which is connected, non-bipartite,
and which is a Ramanujan graph.

   Or $p$ is not a square modulo $q$.
Then $S_{p,q} \cap PSL_2(q) = \emptyset$
and $S_{p,q}$ generates $PGL_2(q)$.
By definition, $X^{p,q}$ is the Cayley graph
of $PGL_2(q)$ with respect to $S_{p,q}$.
It is a  $(p+1)$-regular bipartite graph
with $q(q^2-1)$ vertices
which is connected
and which is a Ramanujan graph.

   See [DaSaV--03] for proofs of
a large part of the facts stated above,
including the connectedness of the graphs
$X^{p,q}$ when $p \ge 5$ and $q > p^8$,
and the expanding property of this family.
For the proof that $\left(X^{p,q}\right)_{q}$
is actually a family
\footnote{
The family is indexed by the set of all odd primes $q$,
and $p$ is a fixed arbitrary odd prime.
}
of Ramanujan graphs,
see the original papers
([LuPhS--88], with a large part
obtained independently in [Margu--88]),
as well as [Sarna--90].

Table I shows the spectrum of $X^{3,q}$ for $q \in \{5,7,11\}$ and
Table II that of $X^{5,q}$ for $q \in \{7,11\}$.
Numerical computations of eigenvalues reported in this paper
have been computed with Mathlab.

\medskip

\begin{prop}
If the graph $X^{p,q}$ is bipartite,
$X^{p,q}$ and its bipartite complement $\overline{X^{p,q}}$
have perfect matchings.
\end{prop}

\emph{Proof}
This is a case of the \lq\lq Marriage Theorem\rq\rq ;
see for example Corollary 1.1.4 in [LovPl--86].
   Here is another reason for $X^{p,q}$ (bipartite {\it or not}):
any connected vertex-transitive graph of even
order has a perfect matching (Section 3.5 in [GodRo--01]);
this applies in particular to Cayley graphs of finite groups of
even order, such as $PGL_2(q)$ and $PSL_2(q)$.
$\square$

\bigskip

\begin{prop} Let $X = (V,E)$ be a finite graph with $n$
vertices and with eigenvalues
$\lambda_0 \ge \lambda_1 \ge \hdots \ge \lambda_{n-1}$.
Let $F$ be a matching of $X$
[respectively of the complement $\overline{X}$]
and let
$\mu_0 \ge \mu_1 \hdots \ge \mu_{n-1}$
be the eigenvalues of $X-F$ [respectively $X+F$].
Then $|\mu_j - \lambda_j| \le 1$
for $j \in \{0,1,\hdots,n-1\}$.
\end{prop}

\emph{Proof}
Outside diagonal entries,
the adjacency matrix $A_F$ of $(V,F)$
coincides with a matrix of permutation
(the permutation being a non-empty product of transpositions with
disjoint supports, one transposition for each edge in $F$).
Thus $\|A_F\| \le 1$. Here, the norm of a matrix acting on the
Euclidean space $\Bbb R^V$  is the operator norm
$\|A_F\| = \sup\left\{ \|Af\|_2 \, \mid \, f \in \Bbb R^V, \
\|f\|_2 \le 1 \right\}$,
where $\|f\|_2^2 = \sum_{v \in V} |f(v)|^2$.

Thus Proposition 1 follows from the classical
Courant-Fischer-Weyl minimax principle,
according to which eigenvalues of  symmetric operators
are norms of appropriate restrictions of these operators.
See e.g. Chapter III in [Bhati--97].
$\square$

\bigskip

\section{Tables}

\medskip
There are several standard efficient algorithms
to find a perfect matching $F$ in a graph $X$;
see [LovPl--86] and [West--01], among others.
We will not describe here the details
of the algorithm we have used.
Eigenvalues of $X-F$ can then
be computed with Mathlab.

   The eigenvalues of a graph of the form $X^{p,q}-F$
depend on the choice of $F$.
Table III gives for each of three pairs $(p,q)$
the values of the spectral gaps $p-\lambda_1(X^{p,q}-F)$
corresponding to four different $F$.
Table III shows that there are situations ($p = 5, q = 7$) with
$\lambda_0(X-F) = k -1 < \lambda_0(X) = k$ and
$\lambda_1(X-F) > \lambda_1(X)$.

Table IV shows the full spectrum of $X^{3,5}-F$ for one specific $F$.
Tables V to VII show the ten largest eigenvalues of three graphs of the
form $X^{p,q} + F$.
Observe that the multiplicities in Tables IV to VII are much less than
those  of the unperturbed graphs. 
\newline\newline\newline\newline\newline\newline\newline\newline\newline
\newline\newline\newline\newline\newline\newline\newline\newline

\vskip.5cm
\begin{center}

\newpage

\begin{tabular}{|c|c||c|c||c|c|}
\hline
\multicolumn{6}{|c|}{Table I: spectra of  $X^{3,q}$}\\
\hline
\multicolumn{2}{|c||}{q=5}&\multicolumn{2}{|c||}{q=7}&\multicolumn{2}
{|c|}{q=11}\\
\hline
eigenvalues &multiplicities &eigenvalues &multiplicities &eigenvalues 
&multiplicities\\
\hline
   -4.0000  &  1& -4.0000 &   1&  -3.2361 &  30\\
   -3.0000  & 12&-3.0000 &  24& -3.0000  & 33\\
   -2.0000  & 28&-2.8284  & 30 & -2.7321  & 10\\
   -1.0000  &  4& -2.0000  & 28&  -2.6180 &  24\\
    0.0000  & 30&-1.4142  & 24 &  -2.3723 &  10 \\
    1.0000  &  4& -1.0000  & 40&  -2.0468 &  36  \\
    2.0000  & 28&  0.0000  & 42&  -2.0000 &  10 \\
    3.0000  & 12& 1.0000  & 40 &  -1.6180 &  36   \\
    4.0000  &  1 &1.4142 &  24 &   -1.5616&   33   \\
            &         & 2.0000 &  28&    -0.9191&   36  \\
            &         &    2.8284  & 30&  -0.7321&   30   \\
   &      & 3.0000   &24               &  -0.3820&   24  \\
 & &   4.0000    &1                      &   0.0000&   30   \\
    &       &                     &           &0.3820 &  12\\
                        & &                &  &0.6180  & 36\\
                               &          & & &0.7321  & 10\\
                                        & & & & 1.0000 &  52\\
                                       & &  & &  1.2361&   30\\
                                       & &  & &1.9191  & 36\\
                                        & & & &2.0000  & 20\\
                                        & & & &2.5616  & 33\\
                                        & & & & 2.6180 &  12\\
                                        & & & &2.7321  & 30\\
                                        & & & & 3.0468 &  36\\
                                        & & & &3.3723  & 10\\
                                        & & & & 4.0000 &   1\\

\hline
\end{tabular}
\end{center}

\begin{center}
\begin{tabular}{|c|c||c|c|}
\hline
\multicolumn{4}{|c|}{Table II: spectra of $X^{5,q}$}\\
\hline
\multicolumn{2}{|c||}{q=7}&\multicolumn{2}{|c|}{q=11}\\
\hline
eigenvalues &multiplicities & eigenvalues &multiplicities\\
\hline

  -6.0000  &  1& -4.0243&   36\\
   -4.0000  & 21& -3.7321&   30\\
   -3.0000   &16& -3.0000&   65\\
   -2.8284   &42& -2.2361&   30\\
   -2.0000&   21& -1.7321&   10\\
   -1.4142 &  12& -1.6180&   60\\
    -1.0000 & 48& -1.3723&   10\\
    0.0000   &14& -1.2361&   12\\
    1.0000   &48& -0.5616&   33\\
    1.4142   &12& -0.2679&   30\\
    2.0000   &21& -0.1638&   36\\
     2.8284  &42& 0.6180 &  60\\
     3.0000  &16& 1.0000 &  30\\
    4.0000   &21& 1.7321 &  10\\
    6.0000  &  1&  1.7818&   36\\
             &      &   2.2361&   30\\
              &      &  3.0000&   50\\
              &      &  3.2361&   12\\
              &      &  3.4063&   36\\
              &      &  3.5616&   33\\
              &      &  4.3723&   10\\
              &      &  6.0000  &  1\\

\hline
\end{tabular}
\end{center}

\begin{center}
\begin{tabular}{|c|c|c|}
\hline
\multicolumn{3}{|c|}{Table III: spectral gaps for $X^{p,q}-F$}\\
\hline
\multicolumn{1}{|c|}{p=3,q=5}&\multicolumn{1}{|c|}{p=3,q=7}&
\multicolumn{1}{|c|}{p=5,q=7}\\
\hline
0.4457&   0.2499&  0.7910\\  
0.3025 & 0.1862&   0.7732\\
0.2993 &  0.1785 &   0.7367\\ 
0.2702&     0.0272&   0.7152  \\
\hline
\end{tabular}
\end{center}

\begin{center}
\begin{tabular}{|c|c|}
\hline
\multicolumn{2}{|c|}{Table IV: spectrum of $X^{3,5}-F$}\\
\hline
eigenvalues &multiplicities \\
\hline

    -3.0000 &   1\\
   -2.5543  &  8\\
   -2.5450  &  4\\
   -2.1542  &  4\\
   -2.0000  &  6\\
   -1.8829  &  8\\
   -1.2929  &  8\\
   -1.0000  &  3\\
   -0.8302  &  4\\
   -0.5086  &  8\\
   -0.4394  &  4\\
    0.0000  &  4\\
    0.4394  &  4\\
    0.5086  &  8\\
    0.8302  &  4\\
    1.0000  &  3\\
    1.2929  &  8\\
    1.8829  &  8\\
    2.0000  &  6\\
    2.1542  &  4\\
    2.5450  &  4\\
    2.5543  &  8\\
    3.0000  &  1\\
\hline
\end{tabular}
\end{center}

\begin{tabular}{|c|c||c|c||c|c|}
\hline
\multicolumn{6}{|c|}{Table V: largest eigenvalues for
$X^{3,5}+F$}\\
\hline
eigenvalues &multiplicities &eigenvalues &multiplicities &eigenvalues
&multiplicities\\
\hline
  3.2578    & 1 & 3.2163&1&3.1707 &1\\
    3.3225  &  1&3.3208 &1& 3.1998 &1 \\
    3.3425  &  1& 3.3431 &1&    3.2214&1\\
    3.4295  &  1&3.4417&1&3.2418  &1\\
    3.4859  &  1& 3.4992&1&  3.3046&1 \\
    3.5140  &  1& 3.5358& 1&    3.5525&1\\
    3.5687  &   1&3.6211 & 1&  3.5653 &1\\
    3.5950  &  1& 3.6822  & 1& 3.5935 &1\\
    3.6758  &  1& 3.8466   &1& 3.6547  &1 \\
    5.0000  &  1 &5.0000  &1&  5.0000  &1\\

\hline
\end{tabular}
\newpage
\begin{tabular}{|c|c||c|c||c|c|}
\hline
\multicolumn{6}{|c|}{ Table VI: largest eigenvalues for
$X^{3,7}+F$}\\
\hline
eigenvalues &multiplicities &eigenvalues &multiplicities &eigenvalues
&multiplicities\\
\hline

    3.6042&    1& 3.6199  &1& 3.6138  &1\\
    3.6130 &   1& 3.6478 &1 &3.6431 &1\\
    3.6349  &  1&3.6594   &1&  3.6524&1 \\
    3.6728   & 1&3.6826   &1& 3.6726  &1\\
    3.6892    &1& 3.6996  &1&3.6922   &1\\
    3.6971&    1&3.7203  &1& 3.7131   &1\\
    3.7073 &   1&3.7468   &1& 3.7275 &1\\
    3.7505  &  1& 3.7548  &1& 3.7461  &1\\
    3.7697   & 1&3.7752 &1& 3.7985    &1\\
    5.0000   &1& 5.0000 &1&  5.0000&1\\

\hline
\end{tabular}

 \vspace{1cm}
 \begin{tabular}{|c|c||c|c||c|c|}
\hline
\multicolumn{6}{|c|}{ Table VII: largest eigenvalues for
$X^{5,7}+F$}\\
\hline
eigenvalues &multiplicities &eigenvalues &multiplicities &eigenvalues
&multiplicities\\
\hline
4.3702    &1& 4.3388&1& 4.3229    &1\\
    4.4015 &   1& 4.3738&1&  4.3405&1\\
    4.4271  &  1& 4.4326 &1& 4.3882 &1 \\
    4.4625   & 1& 4.4790  &1&  4.4117&1 \\
    4.4888    &1& 4.5124  &1& 4.4671 &1\\
    4.4971&    1& 4.5618  &1& 4.5585&1\\
    4.5819&    1&  4.5925  &1&4.5875 &1 \\
    4.5976 &   1& 4.6417  &1& 4.6341&1\\
    4.6512  &  1& 4.6892   &1& 4.7260&1\\
    7.0000&1&  7.0000  &1& 7.0000 &1\\

\hline
\end{tabular}

\bigskip
{\huge References}
\medskip
\begin{description}
\item[[Bhati--97]]  R. Bhatia,
\emph{ Matrix analysis},
 Graduate Texts in Mathematics 169, Springer  1997

\item[[Colin--98]]  Y. Colin de Verdi\`ere,
\emph{ Spectres de graphes},
 Cours sp\'ecialis\'es {\bf 4}, Soc. Math. France  1998.

\item[[DaSaV--03]] G. Davidoff, P. Sarnak, and A. Valette,
\emph{ Elementary number theory, group theory, and Ramanujan graphs}, London Math. Soc. Student Texts {\bf 55}, Cambridge Univ. Press 
2003.

\item[[Fried--91]] J. Friedman, 
\emph{On the second eigenvalue and random walks in random 
$d$-regular graphs}, Combinatorica 11 (1991) 331--362.

\item[[GodRo--01]]  C. Godsil and G. Royle,
\emph{ Algebraic graph theory},
 Graduate Texts in Mathematics 207, Springer  2001.

\item[[LovPl--86]]  L. Lov\'asz and M.D. Plummer,
\emph{ Matching theory},
 Annals Discrete Math. {\bf 29}, North Holland,
 1986.

\item[[Lubot--94]]  A. Lubotzky,
\emph{ Discrete groups, expanding graphs and invariant measure},
 Birkh\"auser  1994.

\item[[LuPhS--88]]  A. Lubotzky, R. Phillips and P. Sarnak,
\emph{ Ramanujan graphs},
Combinatorica {\bf 8},  1988, pages 261--277.

\item[[Margu--88] ]  G.A. Margulis,
\emph{ Explicit group-theoretical constructions of combinatorial
schemes and their application to the design of expanders and
concentrators},
J. Probl. Inf. Transm {\bf 24}, 1988, pages 39--46.

\item[[Sarna--90]]  P. Sarnak,
\emph{ Some applications of modular forms},
 Cambridge University Press  1990.

\item[[Valet--97]]  A. Valette,
\emph{Graphes de Ramanujan et applications},
 S\'eminaire Bourbaki, expos\'e 829, Ast\'e\-risque,
 {\bf 245}, Soc. Math. France 1997, pages 247--296.
 
\item[[West--01]] D.B. West,
{\it Introduction to graph theory}, second edition,
Prentice Hall 2001.

\end{description}

\end{document}